\documentclass{article}

\usepackage{amssymb}

\usepackage[T1]{fontenc}
\usepackage[dvips]{graphicx}

\newtheorem{theorem}{Theorem}[section]
\newtheorem{definition}[theorem]{Definition}
\newtheorem{lemma}[theorem]{Lemma}
\newtheorem{corollary}[theorem]{Corollary}
\newtheorem{proposition}[theorem]{Proposition}
\newtheorem{conjecture}[theorem]{Conjecture}

\newcommand{\PROOF}{\noindent {\bf Proof}: }

\begin{document}

\title{ON BLOCKING NUMBERS OF SURFACES}
\author{Wing Kai, Ho}
\date{July 18, 2008}

\maketitle

\begin{abstract}
The blocking number of a manifold is the minimal number of points
needed to block out lights between any two given points in the
manifold. It has been conjectured that if the blocking number of a
manifold is finite, then the manifold must be flat. In this paper we
prove that this is true for 2-dimensional manifolds with non-trivial
fundamental groups.
\end{abstract}

\section {Introduction}
It has been asked whether how lights on a manifold being blocked
determines the geometry of the manifold. In this paper we prove a
conjecture about blocking numbers for the 2-dimensional case.\\

Throughout the whole paper, we let $(M,g)$ to be a smooth, closed,
orientable 2-dimensional Riemannian manifold, where $g$ is the
Riemannian metric. By a \emph {geodesic segment} we mean a geodesic
$\gamma: [0,a] \rightarrow M$ where $a$ is the length of
$\gamma$.\\

Given two points on $(M,g)$, we can connect them by geodesic
segments. We define the \emph {blocking number} between the two
points to be the minimal number of points needed to block all
geodesic segments connecting them. The blocking number of the
manifold $(M,g)$ is then defined to be the supremum of the blocking
numbers between any two points in $(M,g)$. In some literatures, this
number is referred as the \emph {security threshold} of $(M,g)$. A
manifold
with finite blocking number is said to be \emph {secure}.\\

It is suggested that the blocking number of a Riemannian manifold
could give some information about the geometry of that manifold. A
general conjecture is that if $(M,g)$ has finite blocking number,
then the metric $g$ must be flat. Our task in this paper is to show
that this conjecture is true for orientable surfaces with
non-trivial fundamental groups.\\

Similar results for 2-dimensional tori have been obtained
independently by V. Bangert and E. Gutkin in a preprint posted on
ArXiv on 22 Jun 08, using an analogous but a somewhat different
method, please see \cite{BaG}.\\

\section {Preliminaries and Previous Results}
Let us assume that $(M,g)$ is a smooth, complete, compact orientable
Riemannian manifold. We also assume that geodesics are parametrized
by arc length. If $\gamma: [s_0,s_1] \rightarrow M$ is a geodesic
segment, we call $\gamma(s_0)$ and $\gamma(s_1)$ the \emph
{endpoints} of $\gamma$ and all other points of $\gamma$ the \emph
{interior points} of $\gamma$. We say a geodesic segment \emph
{connects} two points $x$ and $y$ if the points are the endpoints of
$\gamma$. A geodesic $\gamma$ is said to be \emph {blocked} by a
point $z$ if $z$ is an interior point of $\gamma$.\\

Given any two points in $(M,g)$, we can define the blocking number
between them as follows:\\

\begin{definition}
Let $x$ and $y$ be two points in $(M,g)$. The blocking number
$B(x,y)$ between $x$ and $y$ is a positive integer (could be
infinite) that is the minimal number of points needed to block away
all geodesic segments connecting $x$
and $y$.\\
\end{definition}
Next we can give the definition of the blocking number of a
manifold.\\

\begin{definition}
The blocking number of a Riemannian manifold $(M,g)$, denoted by
$B(M,g)$, is the supremum of the blocking numbers between any two
points in $M$. i.e.,\\

\centerline{$B(M,g)=\sup \{B(x,y)|~x,y\in M\}$}
\end{definition}

\vspace{.15in}

 {\bf Example 1}: Let $(M,g)$ be a Hadamard manifold,
which means $(M,g)$ is simply connected and has non-positive
sectional curvature everywhere, this way we have $B(M,g)=1$. This is
because given any two points in $(M,g)$, the Cartan-Hadamard theorem
implies that these two points
are connected by a unique geodesic segment.\\

{\bf Example 2}: Let $(M,g)$ be a flat torus, then $B(M,g)=4$. To
see why this is the case, let us first remark that blocking number
is invariant under affine transformation, therefore we can assume
that $(M,g)$ is the standard flat torus given by $\mathbb R^2
/\mathbb Z^2$. Now we want to show that for any two points $x,y\in
\mathbb R^2 /\mathbb Z^2$ we have $B(x,y)=4$. Without loss of
generality, let us assume $x$ be a point in $M$ and $\tilde{x} \in
\mathbb R^2$ with $\tilde{x}=(0,0)$, such that $\tilde{x}$ projects
to $x$. Now if $y$ is any other point in $M$ with coordinates
$(a,b)\in [0,1]\times [0,1]$, then $y$ is lifted to points in
$\mathbb R^2$ with coordinates $(a+m,b+n)$, $m,n\in \mathbb Z$. Let
us connect each $(a+m,b+n)$ to $\tilde{x}$ by a straight segment
$\tilde{\gamma}_{m,n}$. It is easy to see that the projections of
$\tilde{\gamma}_{m,n}$ to $M$ coincide with all geodesics connecting
$x$ and $y$. Next, for each $(m,n)$ let us consider the midpoint of
$\tilde{\gamma}_{m,n}$ to be its the blocking point. These blocking
points projects to $M$ having coordinates
$(\frac{a}{2},\frac{b}{2})+(\frac{1}{2} \mathbb Z)^2/\mathbb Z^2$.
It is easy to see there could be at most 4 of such blocking points.
Therefore we conclude that the set of geodesics connecting $x$ and
$y$ can be blocked by 4 points, since $x$ and $y$ are
arbitrary, the blocking number of $M$ is 4.\\

{\bf Example 3}: Let $(M,g)$ be a standard $n-sphere$, then
$B(M,g)=\mathfrak{c}$, where $\mathfrak{c}$ stands for the
cardinality of the continuum. This is because if $x$ and $y$ are two
non-antipodal points, then they are connected by exactly two
geodesics, such that these geodesics form a great circle. Now if $x$
and $y$ are two antipodal points, then they are connected by a
family of distinct geodesics, each of which is a half great-circle.
This family of geodesics has cardinality of $\mathfrak{c}$ and so
$B(M,g)=\mathfrak{c}$.\\

The consideration of blocking number could originate from the study
of polygonal billiard systems and geometric optics. One interesting
question is that how the geometry of the manifold could relate to
its blocking number. To begin with the discussions, let us remark a
result that states flat manifolds have finite blocking
numbers \cite{BS}:\\

\begin{theorem}
Compact flat manifolds have finite blocking numbers.
\end{theorem}

A natural question to ask is if the converse of the above theorem is
true.\\

\begin{conjecture}
A compact Riemannian manifold has finite blocking number if and only
if it is flat.
\end{conjecture}

In \cite{BG}, K. Burns and E. Gutkin have obtained a partial
solution to the conjecture. They have related the blocking
properties with metric entropy of the geodesic flow.\\

\begin{theorem}
If $M$ is a manifold without conjugate points, and the geodesic flow
of $M$ has positive metric entropy, then $B(M)$ is not finite.
\end{theorem}

Using this theorem, we can obtain the following result for manifolds
with everywhere non-positive curvature \cite{LS}.\\

\begin{corollary}
Let $(M,g)$ be a Riemannian manifold of non-positive curvature. If
the blocking number of $(M,g)$ is finite, then $(M,g)$ is flat.
\end{corollary}

\PROOF Assume that $(M,g)$ has non-positive curvature and is not
flat, then $(M,g)$ has no conjugate or focal points. According to
Pesin \cite{Pe}, this means the geodesic flow of $(M,g)$ has
positive entropy, hence by Theorem 2.5, the blocking number of $M$
is not finite.$\square$\\

However the assumption of having no conjugate points is rather
strong. For instance if $M$ is a torus and we assume that $M$ has no
conjugate point, this automatically means that $M$ is flat by the
Hopf's conjecture \cite{BI} and the blocking number cannot play a
role
here.\\

\section{The Main Result}

We are now ready to proceed to the main theme of this paper. Our
main result is that Conjecture 2.4 is true for two dimensional
surfaces with non-trivial fundamental groups:\\

\begin{theorem}
Let $(M,g)$ be a compact, complete, orientable 2-dimensional
Riemannian manifold with $\pi_1(M) \neq \{0\}$, then $(M,g)$ has
finite blocking number if and only if $(M,g)$ is a flat torus.
\end{theorem}

In view of Theorem 2.3, we only need to prove the only if part of
the theorem. First of all, let us recall the following
classification theorem for compact orientable 2-dimensional
Riemannian manifolds.\\

\begin{theorem}
The homeomorphism classes of 2-dimensional manifolds are determined
by the genus.
\end{theorem}

Therefore what we will do is to separate the surfaces according to
their genus. In the following we will prove some key lemmas and then
we will use the lemmas to investigate surfaces of different
genus.\\

\section{Key Tools and Lemmas}

We now provide all notations and lemmas needed to prove Theorem 3.1.
Let us start by recalling some fundamental concepts in Riemannian
geometry.\\

Recall that $(M,g)$ is a geodesic complete, locally compact
Riemannian manifold. If we denote by $L_g(\gamma)$ the length of the
curve $\gamma$, We can define the following metric
$d_g(\cdot,\cdot)$ on $M$:\\

\centerline{$d_g(x,y)=\inf_{\gamma} \{\gamma:[a,b]\rightarrow M$
continuous, $\gamma(a)=x,\gamma(b)=y\}$}

\vspace{.15in}

In this way $(M,d_g)$ can be considered a complete metric space. A
continuous curve in $M$ is said to be \emph{distance minimizing} if
the distance in terms of $d_g$ between the end points is equal to
the length of the curve. Next, we define a \emph{minimal geodesic}
to be a continuous curve that is distance minimizing on its
subsegments:\\

\begin{definition} A minimal geodesic is a geodesic segment
$\gamma:[a,b] \rightarrow M$ such that for all $s,s'\in [a,b]$ we
have,\\

\centerline{$d_g(\gamma(s),\gamma(s'))=|s-s'|$}

\end{definition}
In particular, the Arzela-Ascoli theorem, completeness and
compactness of $(M,d_g)$ gives us the following version of the
Hopf-Rinow theorem:\\

\begin{theorem}
Any two points of $M$ can be connected by a minimal geodesic.
\end{theorem}

It is easy to see that $d_g$ is realized by the length of the
shortest minimal geodesic, therefore we can call $d_g$ the
\emph{geodesic distance} on $M$.\\

Now if we consider closed curves, i.e. $\gamma:[a,b]\rightarrow M$
satisfying $\gamma(a)=\gamma(b)$. A closed curve that is a geodesic
is called closed geodesic or periodic geodesic. A periodic geodesic
is minimal if it unwraps to a minimal geodesic.\\

\begin{definition}
Let $\gamma$ be a periodic minimal geodesic and $x$ be a point not
on $\gamma$, then the \emph{distance} between $x$ and $\gamma$,
denoted by $d(x,\gamma)$, is the shortest geodesic distance between
them as subsets, i.e.\\

\centerline{$d(x,\gamma)=\inf \{d_g(x,y)|y\in \gamma\}$}
\end{definition}

\begin{definition}
Let $\gamma$ be a periodic minimal geodesic and $c$ be a minimal
geodesic. $c$ is said to be \emph{asymptotic} to $\gamma$ if $c$
does not intersect $\gamma$ and for all $\epsilon>0$, there exists
$s_0>0$ such that\\

\centerline{$d(c(s),\gamma) <\epsilon$, whenever $s>s_0$}

\end{definition}

The first lemma is a simple fact in differential geometry:\\

\begin{lemma}
Two minimal geodesics originating from the same point will not be
minimal beyond their first point of intersection.
\end{lemma}

For a proof of the lemma, see for instance, \cite{CE}.\\

We want to relate the blocking number of $M$ with periodic minimal
geodesics. Next we will present a crucial lemma that shows the
condition such that $M$ must have infinite blocking number. Let us
first remark that if $\gamma$ is a non-trivial periodic minimal
geodesic,
then for any sufficiently small $\epsilon>0$, the set\\

\centerline{$U_\epsilon(\gamma):=\{x\in M|d(x,\gamma) \leq
\epsilon\}$}

\vspace{.15in}

\noindent has $\gamma$ as its retract. Note that
$U_\epsilon(\gamma)$ is diffeomorphic to an annulus. So if we
consider the universal cover $\tilde{M}$ of $M$,
$U_\epsilon(\gamma)$ is lifted to a infinite strip. Also $\gamma$
will be lifted to an infinitely long minimal
geodesic $\tilde{\gamma}$.\\

\begin{lemma}
If M has a closed minimal geodesic and another minimal geodesic
asymptotic to it, then M has infinite blocking number.
\end{lemma}

\PROOF Suppose the contrary is true, that if $x,y$ are two points in
$M$, we have $B(x,y)<\infty$. We show that this could lead to a
contradiction by demonstrating that there is a point $x$ on the
asymptotic minimal geodesic $c$ stated above and a point $y$ on the
closed geodesic $\gamma$, such that the blocking number between $x$
and $y$ cannot be finite. To prove this claim, let $y$ be any point
on $\gamma$, then consider the lifts $\tilde{\gamma}$ of $\gamma$.
$y$ would be lifted to a countable set of points $y_i$ in a long
strip. For any given point $x$ on $c$ we can connect $x$ to $y_i$ by
a minimal geodesic $\tilde{\gamma_i}$. Each of these
$\tilde{\gamma_i}$ projects to $M$ a geodesic $\gamma_i$
connecting $x$ and $y$.\\

If the point $x$ is sufficiently close to $\gamma$, then there must
be infinitely many $\gamma_i$ such that each of them wraps around
$\gamma$ in the same direction as the asymptotic geodesic $c$. We
see that no two of these $\gamma_i' s$ can intersect at points other
than $x$ and $y$. This is because if $\gamma_i$ and $\gamma_j$ do
intersect, then their lifts $\tilde{\gamma_i}$ and
$\tilde{\gamma_j}$ will intersect before they hit $y_i$ and $y_j$
respectively. This contradicts the fact that two minimal geodesics
originating from the same point will not be minimal beyond their
first point of intersection. For the same reason, each of these
geodesic $\gamma_i$ can only intersect $c$ at $x$.\\

So let us fix $x$ and $y$ as above and show that $B(x,y)=\infty$. If
not, let us assume $\{z_i\}$ to be a finite set of points that block
all geodesics $\{\gamma_i\}$. Throw away several points from
$\{z_i\}$ if necessary, we could assume that all $z_i$ are not on
$\gamma$. Hence there exists a neighborhood $U$ of $\gamma$ in $M$
such that no $z_i$ belongs to $U$. As mentioned above the geodesic
$c$ enters and stay in $U$ after a finite length, and since each
$\gamma_i$ cannot intersect $c$ before it hits the point $y$, there
is a uniform $t_0$ such that for each $i$, $\gamma_i(t)$ is in $U$
for all $t>t_0$.\\

So the infinite set of geodesics $\gamma_i|_{[0,t_0]}$, which only
intersect at $x$, have to be blocked by a finite set of points
$\{z_i\}$. This is a contradiction and the lemma is
proven.$\square$\\

Our next tool is monotone twist maps, which can be used to analyze
geodesic behaviors on a 2-dimensional torus, details of monotone
twist maps will be covered in the next section. For now let us see
how the previous lemma can prove a weaker result concerning the
blocking properties of a 2-dimensional torus.\\

Let us now assume that $(M,g)$ is a 2-dimensional torus. The metric
$g$ is said to be \emph{bumpy} if all periodic geodesics of $(M,g)$
are non-degenerate. We will now show that a bumpy torus has infinite
blocking number.\\

\begin{proposition}
If $(M,g)$ is a bumpy torus, then it has infinite blocking number.
\end{proposition}

\PROOF Firstly, let $\gamma$ be a shortest periodic geodesic of $M$.
Then $\gamma$ is a hyperbolic geodesic due to Morse \cite{Mo}. This
implies that $\gamma$ is isolated. Now we pick another periodic
minimizing geodesic $\gamma'$ in the same homotopy class, such that
there is no periodic geodesic from this homotopy class lies in the
annulus bounded by $\gamma$ and $\gamma'$. This is possible since
$\gamma$ is isolated, if $\gamma$ is the only closed geodesic in the
homotopy class that we can consider $\gamma'=\gamma$ such that the
lift of $\gamma'$ is next to the lift of $\gamma$. Now we observe
that the annulus bounded by $\gamma$ and $\gamma'$ does not contain
any periodic minimal geodesics since geodesics from all other
homotopy class must either intersect $\gamma$ or itself, but we know
that a geodesic that intersect itself cannot be minimizing.
Therefore we conclude that for the strip in $\mathbb R^2$ between
the lift of $\gamma$ and $\gamma'$, there cannot be a curve which
is the lift of a periodic minimal geodesic.\\

Using monotone twist maps, Bangert \cite{Ba} [Th 6.8] proved that
there exists a minimal geodesic such that $c$ is $\omega$-asymptotic
to $\gamma$. So in particular, if $U$ is a neighborhood of $\gamma$
in $T$, $c$ would stay inside $U$ after a finite length. By Lemma
4.6 we can conclude that $(M,g)$ cannot have finite blocking
number.$\square$\\

Note that Proposition 4.7 is an immediate corollary of our main
result.\\

\section{Proof of The Main Result}

We now furnish the proof of Theorem 3.1. In view of Theorem 3.2, we
could separate the surfaces in term of their genus $\mathfrak{g}$.
We will first prove that when $\mathfrak{g}=1$, the only case that
$M$ has finite blocking number is the flat torus. We will then prove
that when $\mathfrak{g}>1$, $M$
cannot have finite blocking number.\\

\subsection{genus $\mathfrak{g}=1$}

Let $(M,g)$ be a topological 2-dimensional torus, we want to show
that if $(M,g)$ has finite blocking number, than $g$ must be a
flat metric.\\

Our approach is the following: We will argue by first assuming that
$(M,g)$ has finite blocking number, next for each free homotopy
class $[\alpha]$ of $M$, we call an annulus \emph{bad} if the
annulus is bounded by two periodic minimal geodesics from $[\alpha]$
such that no other periodic geodesics exists in the annulus. We will
show that no bad annulus exists in $(M,g)$. Afterwards we prove that
this means $M$ can be foliated by periodic minimal geodesics of the
class $[\alpha]$. This in turns will imply $(M,g)$ has no conjugate
points
and so by Hopf's Theorem $g$ must be a flat metric.\\

Let us begin by introducing some notations. If $\gamma(s)$ is a
periodic minimal geodesic in $M$, then its lift $\tilde{\gamma}(s)$
is a minimal geodesic in $(\mathbb R^2,\tilde {g})$. Let us write
$\tilde{\gamma}(s)=(\xi(s),\eta(s))$ where $\xi$ and $\eta$
represents the coordinates of $\tilde{\gamma}$ in $\mathbb R^2$. We
can define the \emph{average slope} $\alpha$ of $\tilde{\gamma}$ by
setting $\alpha(\tilde{\gamma}):=\lim_{|s| \rightarrow \infty}
\eta(s)/\xi(s)$.\\

Here is a property of the functional $\alpha$ proven in \cite{Ba}:\\

1) If $\xi(s)$ is surjective then $\alpha(\tilde{\gamma})$ exists in
$(-\infty,\infty)$. If $\xi(s)$ is not surjective then $\xi(s)$ is
bounded and we will define $\alpha(\tilde{\gamma})=\infty$.\\

\noindent {\bf Definition}: Let $(q,p)\in \mathbb Z^2$, denote by
$T_{(q,p)}(x)$ the action of the group
$\mathbb Z^2$ on $\mathbb R^2$ that translate the point $x$ by $(q,p)$.\\

\noindent {\bf Definition}: A lifted periodic geodesic
$\tilde{\gamma}$ has \emph{period} $(q,p)\in \mathbb Z^2-\{0\}$ if
the translation $T_{(q,p)} \tilde{\gamma}$ and $\tilde{\gamma}$
coincide up to parametrization. The \emph{minimal period} of a
geodesic is the pair $(q,p)$ such that the geodesic has period
$(q,p)$ and $p,q$ are relatively prime. Obviously, if a periodic
minimal geodesic $\tilde{\gamma}$ has period $(q,p)$, then
$\alpha(\tilde{\gamma})=p/q$.\\

Note that if two periodic minimal geodesics are in the same homotopy
class then their lifts have the same minimal period. It is also easy
to see that two distinct periodic minimal geodesics with the same
period do not intersect [6.6 of \cite{Ba}]. In the following, we
will assume a periodic geodesics $\gamma$ are equivalent to its
iterates $n\gamma, n\in \mathbb N$.\\

Now we are ready to show that a torus with finite blocking number is
flat. Firstly let us state a lemma that is similar to Lemma 4.6.\\

\begin{lemma}
If $(M,g)$ has two periodic minimal geodesics from the same homotopy
class, such that the annulus bounded between does not contain other
periodic minimal geodesic, i.e. if $M$ has a bad annulus, then the
blocking number of $(M,g)$ is infinite.
\end{lemma}

\PROOF Denote one of the periodic minimal geodesics by $\gamma$.
According to 6.8 of \cite{Ba}, there exists a minimal geodesic such
that $c$ is $\omega$-asymptotic to $\gamma$. So in particular, $c$
is asymptotic to $\gamma$ according to Definition 4.4. Therefore we
can apply Lemma 4.6 to conclude that the blocking number of $(M,g)$
is infinite.$\square$\\

Note that in the above lemma, we can replace the phrase 'does not
contain other periodic minimal geodesics' by 'does not contain other
periodic minimal geodesics of the same homotopy class'. It is
because periodic geodesics of other homotopy classes must intersect
either one of the boundaries of the annulus or is a higher
iteration of the geodesic itself.\\

Let $[\alpha] \in \pi_1(M)$, by the first variation formula we know
that there is at least one periodic minimal geodesic. The following
proposition shows that if no periodic minimal geodesic in $[\alpha]$
is isolated, then $M$ is foliated by periodic minimal geodesics from
this homotopy class.\\

\begin{proposition}
Fix $[\alpha] \in \pi_1(M)$, if for any two periodic minimal
geodesics $\gamma$ and $\gamma'$ of $[\alpha]$ there is a periodic
minimal geodesic of $[\alpha]$ that lies in the annulus bounded by
$\gamma$ and $\gamma'$, then $M$ is foliated by periodic geodesics
in $[\alpha]$.
\end{proposition}

\PROOF Let $\gamma$ be the shortest periodic geodesic from the free
homotopy class $[\alpha]$, this closed geodesic lifts to $(\mathbb
R^2, \tilde{g})$ to a minimal geodesic $\tilde{\gamma}$. Let
$\gamma'$ be another periodic minimal geodesic from $[\alpha]$ and
$\tilde{\gamma}'$ be the corresponding lift. If there is no such
geodesic then we let $\gamma=\gamma'$ and $\tilde{\gamma}'$ be the
minimal geodesic in $\mathbb R^2$ neighboring $\tilde{\gamma}$ and
projects to $\gamma$. Next we assert the following is true.\\

{\bf Claim}: For any point $x$ in the strip bounded by
$\tilde{\gamma}$ and $\tilde{\gamma}'$,
there is a minimal geodesic $c_x$ passing through $x$ such that $\alpha(c_x)=\alpha(\gamma)=\alpha(\gamma')$.\\

PROOF of Claim: To prove this claim, assume that the minimal
geodesics $\tilde{\gamma}$ and $\tilde{\gamma}'$ have minimal period
$(q,p)$. Let $x \in \mathbb R^2$ be a point that lies in the strip
between them. Denote by $x_1=T_{(q,p)}x$, and
$x_{i+1}=T_{(q,p)}x_i$, $\forall i\in \mathbb N$. Then we can
connect each $x_i$ to $x$
by one minimal geodesic $c_i$.\\

Note that each $c_i$ cannot touch or cross each of $\tilde{\gamma}$
and $\tilde{\gamma}'$ transversely . It is because if $c_i$ crosses
say, $\tilde{\gamma}$ transversely, then it has to cross it at least
twice, and we know that geodesics that intersect each other twice
cannot be minimal beyond the intersections. If $c_i$ touches
$\tilde{\gamma}$ then this will contradict the
uniqueness of geodesic for a given initial point and tangent vector. This means all $c_i$ stay in the strip.\\

Now, let $v_i$ be a vector in $U_xM$, the unit tangent sphere at
$x$, such that $v_i=c_i'(0)$. Then $\{v_i\}$ is a set of vectors in
the compact sphere. So there is a limit
$v=\lim_{i\rightarrow\infty}v_i$. Let $c(s)$ be the forward complete
geodesic satisfying $c(0)=x$, $c'(0)=v$. We now show $c(s)$ stays
the strip for all $s>0$. Suppose not, say, $c(s)$ intersects
$\tilde{\gamma}$ transversely at some point. Then there is an
$\epsilon>0$ and $s_0>0$ such that $c(s_0)$ is not in the strip and
lay outside of the $\epsilon$-neighborhood of $\tilde{\gamma}$.\\

Since a geodesic is a solution of a second order ODE, the solutions
with initial point $x$ continuously depend on the initial vector
$v$. Also as $M$ is geodesic complete, the geodesic flow
$\phi_s(x,v)$ on the unit tangent bundle of $T$ is defined for all
$s>0$.  Let $f:UM\rightarrow T$ be the composition of the time-$s_0$
map restricted to $U_xM$, $\phi_{s_0}(x,\cdot):U_xM \rightarrow UM$
with the projection $\pi:UM \rightarrow M$ given by $(x,v) \mapsto
x$, i.e. $f=\pi \circ \phi_{s_0}$.  Then $f$ continuously depends on
$v$. So for the geodesic $c$ there exists a $\delta>0$ such that for
all geodesics $\bar{c}$ with $\bar{c}(0)=x$ and
$\bar{c}'(0)=\bar{v}$ for $\|v-\bar{v}\| \leq \delta$, we have
$\tilde{d}(\bar{c}(s_0)-c(s_0)) \leq \epsilon$. Recall that
$v=\lim_{i\rightarrow\infty}v_i$ and $c_i$ is the geodesic at $x$
with initial vector $v_i$, so by the arguments above there exists
arbitrarily long $c_i$ such that $c_i(s_0)$ lies outside the strip.
In particular $c_i$ crosses the boundary of the strip transversely.
However this is not possible since $c_i$ connects two points in the
strip, if it cross the boundary than it will intersect the boundary
minimal geodesic twice, this contradicts the assumption that $c_i$
is a minimal geodesic.\\

This means that $c$ must stay inside the strip, it is now easy to
see that $\alpha(c)=\alpha(\gamma')=\alpha(\gamma)$ and so the claim
true.\\

Now the claim is true so we can then applies Theorem 6.7 of
\cite{Ba}, that $c$ is either periodic or is contained in a strip
between two periodic minimal geodesics $c^-$ and $c^+$. For the
latter case  the strip between $c^-$ and $c^+$ contains no other
periodic minimal geodesic. This will contradict our hypothesis so
$c$ can only be periodic. This means each point $x \in M$ lies on a
periodic geodesic in $[\alpha]$ and the proposition is
proven.$\square$\\

Now we can finalize the proof of Theorem 3.1 for the case when the
genus $\mathfrak{g}=1$.\\

\begin{proposition}
If $(M,g)$ is a 2-dimensional torus with finite blocking number,
then $(M,g)$ is flat.
\end{proposition}

\PROOF Assume that $(M,g)$ has finite blocking number. According to
Lemma 5.1, if $[\alpha]$ is a free homotopy class then between any
two periodic minimal geodesics of $[\alpha]$ there exists another
periodic geodesic. Since the free homotopy class $[\alpha]$ is
arbitrarily, Proposition 5.2 implies that periodic minimal geodesics
in any fixed free homotopy class foliate $M$. So according to Innami
[Corollary 3.2] \cite{In} $(M,g)$ does not have conjugate points and
so by Hopf's Theorem, $g$ must be a flat metric.$\square$\\

\subsection{genus $\mathfrak{g}>1$}

We now let $M$ be a surface of genus $\mathfrak{g}$, where
$\mathfrak{g} \geq 2$. Let $g$ be any metric on $T$. We will now
show that $(M,g)$ must have infinite blocking number:\\

We will make use of Lemma 4.6 again, which states that if a closed
surface has a periodic minimal geodesic and another minimal geodesic
asymptotic to it, then the surface has infinite blocking number.\\

\begin{definition}
Let $\gamma$ be a closed geodesic on $M$. A \emph{collar} of
$\gamma$ is the image of a diffeomorphism $f:\mathbb S^1 \times
[0,a]\rightarrow M$ with $f(\mathbb S^1 \times \{0\})=\gamma$, $a
\in \mathbb R^+$.
\end{definition}
In the following, a \emph{periodic minimal geodesic} is a closed
geodesics with minimal length in the homotopy class. A
\emph{cylinder} is a manifold diffeomorphic to the product of an
interval (open or closed)
and $\mathbb S^1$.\\

We will show that a surface of higher genus has a closed geodesic
and a collar, such that the collar contains an asymptotic geodesic.
Once this is established the proof of our
desired result is trivial.\\

Let $[\alpha]$ be a homotopy class of $M$ which represents the cross
section of a 'handle' of $M$. Using the variation formula we know
there is at least one periodic minimal geodesic $\gamma \in
[\alpha]$. As a matter of fact, periodic minimal geodesics of
$[\alpha]$ are well ordered:\\

\begin{lemma}
Periodic minimal geodesics of $[\alpha]$ form a totally ordered set.
\end{lemma}

\PROOF Recall that no two periodic minimal geodesics from $[\alpha]$
intersect, this is because if they do, then their lifts to the
universal cover will intersect each other more than one time, this
contradicts the assumption that both geodesics are minimal. Since we
also assume that the genus of $M$ is greater than 1, closed
geodesics from $[\alpha]$ cannot homotopically pass the portion of
$M$ containing other handles. Therefore we conclude that periodic
minimal geodesics of $[\alpha]$ form a totally ordered
set.$\square$\\

\begin{proposition}
If M is a surface of genus $g \geq 2$, then there exists a closed
geodesic $\gamma$ and a collar such that except $\gamma$, there is
no closed geodesics homotopic to $\gamma$ that intersect the
collar.\\
\end{proposition}

\PROOF As above, let $[\alpha]$ be a homotopy class of $M$ which
represents the cross section of a 'handle' of $M$. Our claim is that
there is a periodic minimal geodesic of $[\alpha]$ with a collar
such that no geodesics of $[\alpha]$ that intersect the collar.\\

By Lemma 5.4, periodic minimal geodesics of $[\alpha]$ form a
totally ordered set. We also know that periodic minimal geodesics
are the critical points of the length functional, therefore the set
of minimal geodesics is closed. Hence with respect to the ordering
on minimal geodesics of the class $[\alpha]$, there is a 'maximum'
$\gamma$ and 'minimum' $\gamma'$ so that all periodic minimal
geodesics of $[\alpha]$ lies in the cylinder bounded by $\gamma$ and
$\gamma'$.\\

If $\gamma$ and $\gamma'$ are the same geodesic, there could be two
possibilities. The first one is that there is only one periodic
minimal geodesic $\gamma$ in the class $[\alpha]$, in this case we
can find a neighborhood of $\gamma$ such that its retract is
$\gamma$. Since there is only one minimal geodesic, the neighborhood
contains a collar, and the only periodic minimal geodesic contained
in the collar is $\gamma$ itself, thus the statement of the
Proposition is
true in this case.\\

The second possibility is that there are more than one periodic
minimal geodesic. In this case if $\gamma$ and $\gamma'$ are the
same geodesic, then the geodesic $\gamma$ has to homotope along the
handle and overlap with itself. However since the genus of $M$ is
strictly greater than 1, this cannot be done. Therefore we conclude
that for the case where $\gamma$ and $\gamma'$ are the same
geodesic, there could only be one periodic minimal geodesic and we
are done.\\

Now suppose that $\gamma$ and $\gamma'$ are two distinct periodic
minimal geodesics. Note again that $\gamma$ cannot intersect or
touch $\gamma'$. Hence if we denote by $K$ the cylinder bounded by
$\gamma$ and $\gamma'$, there is a neighborhood $U$ of $K$ that is
diffeomorphic to a cylinder. Now since there is no periodic minimal
geodesics from $[\alpha]$ outside of $K$, the set $(U\cup \gamma)
\backslash K$ contains a closed geodesic $\gamma$ and a collar such
that except $\gamma$, there is no closed geodesics homotopic to
$\gamma$ that intersect the collar, finishing the proof of the
Proposition. $\square$\\

Next proposition will direct us to the proof of higher genus case of
Theorem 3.1:\\

\begin{proposition}
If M is a surface of genus $g>1$, then M has a closed minimal
geodesic and another minimal geodesic asymptotic to it.
\end{proposition}

\PROOF According to Proposition 5.5, $M$ has a periodic minimal
geodesic $\gamma$ and a collar $C$ such that except $\gamma$, there
is no closed geodesics homotopic to $\gamma$ that intersects $C$.\\

Note again that $C$ is diffeomorphic to a cylinder. We now cut an
open neighborhood $U$ of $C$ out such that $U$ is diffeomorphic to a
cylinder, then $\gamma$ is the shortest closed geodesic in $U$.
After that we smoothly glue a closed Riemannian cylinder to $U$, we
would then obtain a 2 dimensional Riemannian torus.\\

Since $\gamma$ is the shortest closed curve in the $U$, we can glue
the cylinder to $U$ such that $\gamma$ remains to be the shortest
closed geodesic in the homotopy class. For instance, let the metric
of the cylinder be that all closed curves in it have lengths not
shorter than $\gamma$, which is possible since the length of each
boundary of $U$ is not shorter than $\gamma$. This way we can
guarantee that $\gamma$, as a closed curve of the torus, is the
shortest in its free homotopy class.\\

We know that the geodesic behaviors of a surface is completely
determined by the metric. In particular, the geodesic behaviors of
$C$ as a subset of $M$ is the same as that of $C$ as a subset of the
glued 2 dimensional torus. Now as a subset of the glued torus, $C$
does not contain any periodic minimal geodesics homotopic to
$\gamma$. If we consider the universal cover $\mathbb R^2$, $\gamma$
is lifted to a infinitely long minimal geodesic $\tilde{\gamma}$ and
the collar $C$ is lifted to a strip $\tilde{C}$. Let
$\tilde{\gamma}'$ be the lift of an adjacent periodic minimal
geodesic from the same homotopy class such that $\tilde{\gamma}$ and
$\tilde{\gamma}'$ bound $\tilde{C}$. This could be either one of the
followings:\\

a) $\tilde{\gamma}'$ is the lift of a distinct periodic minimal
geodesic.\\

b)$\tilde{\gamma}'$ is another lift of $\gamma$ which is adjacent to
$\tilde{\gamma}$.\\

In case {\bf a} we claim that $\tilde{\gamma}'$ cannot intersect
$\tilde{C}$. To see why this is true, let us first remark that
periodic minimal geodesics of a torus are exactly the lifts of
closed geodesics which have minimal length in their free homotopy
class, see 6.6 \cite{Ba}.\\

If $\gamma_1$ is a closed curve homotopic to $\gamma$ that lies
completely in $C$, then the length of $\gamma_1$ is strictly greater
than the length of $\gamma$, this is because there is no other
periodic minimal geodesics in $C$.\\

Now we assume that $\gamma_1$ is a closed curve partially lies in
$C$. Note that in the gluing process we ensured that all closed
curve outside $C$ is longer than $\gamma$. Together with the facts
that both boundaries of $C$ is not shorter than $\gamma$, we see
that $\gamma_1$
cannot be shorter than $\gamma$.\\

What matters in both case {\bf a} and {\bf b} is that
$\tilde{\gamma}'$ and $\tilde{\gamma}$ bounds a strip that contains
$\tilde{C}$, and the strip does not contain any periodic minimal
geodesics of the same homotopy class. So we can apply Theorem 6.8 of
\cite{Ba} to show that there is a minimal geodesic $c$ in the strip
which is $\omega$-asymptotic to $\tilde{\gamma}$. Now we consider
the portion of the minimal geodesic $c$ which stays in $\tilde{C}$,
this geodesic projects to a minimal geodesic asymptotic to $\gamma$.
Finally, since the geodesic behaviors of an area is determined
locally by the metric, therefore if we revert the cutting and gluing
procedures that changes $M$ to a 2 dimensional torus, the geodesic
$c$ staying in the collar $C$ remains exactly the same. Therefore on
the surface $M$, $c$ is a minimal geodesic asymptotic to the
periodic minimal geodesic $\gamma$ and the Proposition is
proven.$\square$\\

With the above propositions and lemma, we can now prove the
following proposition:\\

\begin{proposition}
If $M$ is a closed Riemannian surface with genus greater than or
equal to
2, then the blocking number of $M$ is infinite.\\
\end{proposition}

\PROOF By Proposition 5.6, $M$ has a closed minimal geodesic and
another minimal geodesic asymptotic to it. Therefore we can apply
Lemma 4.6 to conclude that $M$ has infinite blocking number.
$\square$\\

Combining the statements of Propositions 5.3 and 5.7, and using
Theorem 3.2, we have proven Theorem 3.1.\\





{\footnotesize \centerline{\rule{9pc}{.01in}}
\bigskip

\centerline{The Pennsylvania State University, University Park, PA
16801} \centerline{e-mail: ho\_wk@math.psu.edu}
\medskip

\end{document}